\documentclass{amsart}
\usepackage[T1]{fontenc}   

\newtheorem{theorem}{Theorem}[section]

\theoremstyle{definition}

\numberwithin{equation}{section}

\begin{document}


\renewcommand{\bf}{\bfseries}
\renewcommand{\sc}{\scshape}
\vspace{0.5in}

\title[non-blockers of singletons, all the possible examples]%
{The hyperspace of non-blockers of singletons, all the possible examples}

\author{Alejandro Illanes}
\address{Instituto de Matem\'aticas, Universidad
Nacional Aut\'onoma de M\'exico, Circuito Exterior, Cd. Universitaria, M\'exico
04510, Ciudad de M\'exico.}
\email{illanes@matem.unam.mx}
\thanks{The first author was supported in part by the projects "Teor\'ia de Continuos,
Hiperespacios y Sistemas Din\'amicos III", (IN 106319) of PAPIIT, DGAPA,
UNAM; and "Teor\'ia de Continuos e Hiperespacios, dos" (AI-S-15492) of CONACYT}

\author{Benjamin Vejnar}
\address{Department of Mathematical Analysis, Faculty of Mathematics and
Physics, Charles University, Czechia}
\email{vejnar@karlin.mff.cuni.cz}
\thanks{The second author was supported by the grant SVV-2020-260583.}

%

\thanks {Orcid numbers, A. Illanes: 0000-0002-7109-4038; B. Vejnar: 0000-0002-2833-5385}

\begin{abstract} Given a metric continuum $X$, a nonempty proper closed subspace $B$ of $X$,
{\it does not block} a point $p\in X\setminus B$ provided that the union
of all subcontinua of $X$ containing $p$ and contained in $X\setminus B$
is a dense subset of $X$. The collection of all nonempty proper closed subspaces
$B$ of $X$ such that $B$ does not block any element of $X\setminus B$ is
denoted by $NB(F_{1}(X))$. In this paper we prove that for each completely
metrizable and separable
space $Z$, there exists a continuum $X$ such that $Z$ is homeomorphic to
$NB(F_{1}(X))$. This answers a series of questions by Camargo, Capul\'in,
Casta\=neda-Alvarado and Maya.
\end{abstract}

\maketitle

\section{\bf Introduction}
A {\it continuum} is a nonempty non-degenerate compact connected metric space.
A {\it subcontinuum} of a continuum $X$ is a nonempty closed connected subspace
of $X$, so one-point sets are subcontinua. For a continuum $X$, we consider
the hyperspaces:

\begin{center}
$2^{X}=\{A\subset X:A$ is a nonempty closed subset of $X\}$,

$C(X)=\{A\in 2^{X}:A$ is connected$\}$, and for each $n\in \mathbb{N}$,

$F_{n}(X)=\{A\in 2^{X}:A$ has at most $n$ points$\}$.
\end{center}

The hyperspace $2^{X}$ is considered with the Hausdorff metric \cite [Definition
2.1]{in}. The {\it hyperspace of singletons} is $F_{1}(X)$.

Given a continuum $X$, a non-empty proper closed subspace $B$ of $X$,
{\it does not block} a point $p\in X\setminus B$ provided that the union
of all subcontinua of $X$ containing $p$ and contained in $X\setminus B$
is a dense subset of $X$. The collection of all nonempty proper closed subspaces
$B$ of $X$ such that for each $z\in X\setminus B$, $B$ does not block $z$
is
denoted by $NB(F_{1}(X))$.

The concept of blockers in continua was introduced in \cite{ik}. Since then,
many aspects of blockers have been studied by several authors in \cite{bpv},
\cite{bpv1}, \cite{cccm}, \cite{cmo}, \cite{elv}, \cite{ik} and \cite{p}.

The particular set of non-blockers $NB(F_{1}(X))$ was studied in \cite{cccm},
where the following general question was considered. 

\bigskip

For what compact metric spaces $Z$ does there exist a continuum $X$ such
that $Z$ is homeomorphic to $NB(F_{1}(X))$?

\bigskip

In \cite{cccm}, this question was answered in the positive for the following
spaces $Z$: a simple closed curve, an arc, each finite set, a convergent
sequence with its limit, the hyperspace $2^{Y}$, for each continuum $Y$,
a $2$-cell, a Cantor set and a Hilbert cube.

The authors of \cite{cccm}, also include the following series of questions.

(a) Are the arc and the simple closed curve the unique $1$-dimensional continua
$Y$ such that $NB(F_{1}(X))$ is homeomorphic to $Y$, for some continuum $X$?

(b) For which finite graphs $Y$ does there exist a continuum $X$ such that
$NB(F_{1}(X))$ is homeomorphic to $Y$?

(c) For which $m\geq 3$ does there exist a continuum $X$ such that
$NB(F_{1}(X))$ is homeomorphic to a simple $m$-od?

(d) For which continua $Y$ does there exist a continuum $X$ such that
$NB(F_{1}(X))$ is homeomorphic to $C(Y)$?

(e) For which continua $Y$ does there exist a continuum $X$ such that
$NB(F_{1}(X))$ is homeomorphic to $Y$?

The aim of this paper is to answer all the questions above by proving the
following.

\begin{theorem}
For each completely metrizable and separable  space $Z$, there is a continuum
$X$ such that $NB(F_{1}(X))$ is homeomorphic to $Z$.
\end{theorem}
\section{\bf The construction}

A mapping is a continuous function.

The letter $P$ will denote the pseudo-arc. For a very complete survey on
the pseudo-arc, see \cite{w}. For properties of the composants of $P$ see
section 2 of Chapter XI of \cite{n1}.     

\bigskip

{\bf Construction of the continuum $X$}

 Let $Q=[0,1]^{\mathbb{N}}$ be the Hilbert cube. Let 
\begin{center}
$Q_{0}=\{x=(x_{1},x_{2},\ldots)\in Q:x_{1}=0\}$.
\end{center} 

We may suppose that $Z$ is a subspace of $Q_{0}$. Since $Z$ is completely
metrizable, it follows that it is a $G_{\delta}$-set in $Q_{0}$. Hence the
set $Q_{0}\setminus Z$ is a countable union of some closed sets $F_{n}\subset
Q_{0}$. Consequently the set $Q\setminus Z$ is the countable union of the
continua: 
\begin{center}
$\{x\in Q:(0,x_{2},x_{3},\ldots)\in F_{n}\}\cup ([2^{-n},1]\times [0,1]^{\mathbb{N}}),
n\in\mathbb{N}$.
\end{center}

Let us denote this collection of continua by $\mathcal{C}$. Hence $Q\setminus
Z=\bigcup \mathcal{C}$.

Fix a countable base $\mathcal{B}$ of $Q$ formed by nonempty open sets whose
closures are continua (in fact, we may suppose that the closures are homeomorphic
to the Hilbert cube).

Fix $p,q\in P$ in different composants.

Let $A$ be the space obtained from $P\times \mathbb {N}$ by identifying the
points $(q,n)$ with $(p,n+1)$, $n\in\mathbb {N}$.

Let $B$ be the space obtained from $P\times \mathbb {Z}$ by identifying the
points $(q,n)$ with $(p,n+1)$, $n\in \mathbb {Z}$.

Clearly $A$ and $B$ are locally compact and $A$ (resp., $B$) has one end
(resp., two ends). Hence for every continuum $K$ (resp., continua $K$ and
$L$) we can compactify $A$ (resp., $B$) in that way that $K$ (resp., disjoint
union of $K$ and $L$) forms the remainder.

We are going to construct the continuum $X$ in such a way that it consists
of the Hilbert cube $Q$ with countably many copies of $A$ and $B$ attached
to it:

$\bullet$ For every $C\in\mathcal {C}$ attach to $Q$ a copy $A_{C}$ of $A$
such that it is limiting to $C$ (i.e. we compactify $A_{C}$ in such a way
that the remainder is $C$) and $A_{C}\cap Q=\emptyset$.

$\bullet$ For every pair $U,V\in\mathcal{B}$, whose closures are disjoint
attach to $Q$ a copy $B_{U,V}$ of $B$ whose one end is limiting to cl$_{Q}(U)$
and the second end to cl$_{Q}(V)$, and $B_{U,V}\cap Q=\emptyset$.

$\bullet$ Arrange this countable attaching process in such a way that $A_{C}$
and $B_{U,V}$ are contained in smaller and smaller neighborhoods of $Q$,
and all the sets in the family 
\[\{A_{C}:C\in\mathcal{C}\}\cup\{B_{U,V}:U,V\in\mathcal{B}, \text{ and } cl_{Q}(U) \cap cl_{Q}(V)=\emptyset\}\] are pairwise disjoint.

The following properties are easy to show.

\begin{enumerate}

\item $X$ is a continuum, 

\item for each $C\in \mathcal{C}$ and for each set $P^{*}=P\times\{n\}$
contained in $A_{C}$, the boundary of $P^{*}$ in $X$ is contained in the
set $\{(p,n),(q,n)\}$. Moreover, the boundary of $A_{C}$ in $X$ is $C$,
and every subcontinuum of $X$ intersecting $A_{C}$ and its complement contains
$C$,

\item for every $U,V\in\mathcal{B}$ with  $cl_{Q}(U)$ $\cap$ $cl_{Q}(V)=\emptyset$
and for each set $P^{*}=P\times\{n\}$
contained in $B_{U,V}$, the boundary of $P^{*}$ in $X$ is the
set $\{(p,n),(q,n)\}$. Moreover, the boundary of $B_{U,V}$ in $X$ is cl$_{Q}(U)$
$\cup cl_{Q}(V)$, and every subcontinuum of $X$ intersecting $B_{U,V}$ and
its complement contains either $cl_{Q}(U)$ or cl$_{Q}(V)$. 

\end{enumerate}

\section{\bf The proof} 

We prove that $NB(F_{1}(X))=F_{1}(Z)$. Take a nonempty closed subset $S$
of $X$ with empty interior. We consider four cases.

{\bf Case 1.} $S$ is not a subset of $Q$. 

In this case there exists a point $s\in S$ which is contained in some copy
$P^{*}=P\times\{n\}$ of $P$ (where $P^{*}$ is contained in some $A_{C}$ or
$B_{U,V}$). Since $S$ has empty interior and $P^{*}\setminus \{(p,n),(q,n)\}$
is open in $X$, there is a point $t\in P^{*}\setminus S$ whose composant
in $P^{*}$ contains neither $(p,n)$ nor $(q,n)$. We claim that $S$ blocks
point $t$. Suppose the contrary, then there exists a subcontinuum $E$ of
$X$ such that $t\in E\subset X\setminus S$ and $E$ is not contained in $P^{*}$. Let $D$
be the component of $E\cap P^{*}$ that contains $t$. By the Boundary Bumping
Theorem \cite[Theorem 20.2]{na}, $\emptyset\neq D$ $\cap$ Bd$_{E}(E\cap P^{*})\subset$
$Bd_{X}(P^{*})\subset \{(p,n),(q,n)\}$. Thus $D$ is a subcontinuum of $P^{*}$
intersecting two composants of $P^{*}$. Hence $D=P^{*}$ and $s\in E\cap S$,
a contradiction. Therefore $S\notin NB(F_{1}(X))$.

{\bf Case 2.} $S\subset Q$ and $S$ contains at least two points $s,t$, $s\neq
t$.

Take sets $U,V\in\mathcal {B}$ with disjoint closures such that $s\in U$ and
$t\in V$. Consider any point $r\in B_{U,V}$. With a similar argument as in
Case 1, it follows that $S$ blocks $r$. Hence $S\notin NB(F_{1}(X))$.

{\bf Case 3.} $S=\{s\}$ and $s\in Q\setminus Z$.

By the construction of $\mathcal {C}$, there is $C\in\mathcal{C}$ with $s\in
C$. Again, as in Case 1, it follows that $S$ blocks every point from $A_{C}$.
Hence $S\notin NB(F_{1}(X))$.

{\bf Case 4.} $S=\{s\}$ and $s\in Z$.

Fix a point $z\in Q\setminus Z$. In order to show that $S$ does not block
any point from $X\setminus S$ is enough to show that for every point
$w$ in $X\setminus S$ there exists a subcontinuum of  $X\setminus S$ containing
$z$ and $w$. In the case that $w\in Q\setminus Z$, take sets $U,V\in\mathcal
{B}$ with disjoint closures such that $w\in U$, $z\in V$ and $s\notin$ $cl_{Q}(U)$
$\cup$ cl$_{Q}(V)$. Then $B_{U,V}$ $\cup$ $cl_{Q}(U)$ $\cup$ $cl_{Q}(V)$ is
a subcontinuum of $X\setminus S$ containing $w$ and $z$. 

In the case that $w\in A_{C}$ for some $C\in \mathcal {C}$, the continuum
$A_{C}\cup C\subset X\setminus S$, contains $w$ and contains a point $w_{1}\in
Q\setminus S$. By the previous paragraph, it is possible to join $w$ with
$z$ by a subcontinuum of $X\setminus S$.

Finally, in the case that $w\in B_{U',V'}$ for some sets $U',V'\in\mathcal{B}$
with disjoint closures. We may suppose that $s\notin $ $cl_{X}(U')$, $w$ belongs
to the copy $P\times \{1\}$ of $P$ in $B_{U',V'}$ and $cl_{X}(V')\subset 
cl_{X}(\bigcup\{P\times\{n\}\subset B_{U',V'}:n<1\})$. Then  the set $cl_{X}(\bigcup\{P\times\{n\}\subset
B_{U',V'}:n\geq1\})$ $\cup$ $cl_{X}(U')$ is a subcontinuum of $X\setminus
S$ and contains $w$ and points in $Q\setminus S$. As in the previous paragraph
we are done.
Therefore $S\in NB(F_{1}(X))$.
This finishes the proof of Theorem 1.1.

\bibliographystyle{plain}

\end{document}